\documentclass[journal]{IEEEtran}
\usepackage{}
\ifCLASSINFOpdf
\else
\fi

\usepackage{graphicx}
\usepackage{dcolumn}
\usepackage{bm}

\usepackage{graphics} 
\usepackage{epsfig} 
\usepackage{times} 
\usepackage{amsmath} 
\usepackage{amssymb}  
\usepackage{accents}
\usepackage{amsfonts}
\usepackage{fancyhdr}
\usepackage{color}
\usepackage{subfigure}

\newtheorem{theorem}{Theorem}
\newtheorem{lemma}[theorem]{Lemma}

\newtheorem{proposition}[theorem]{Proposition}

\newtheorem{definition}{Definition}

\newtheorem{remark}{Remark}

\allowdisplaybreaks[4] 

\begin{document}
%
\title{A {De~Giorgi} Iteration-based Approach for the Establishment of ISS Properties of a Class of Semi-linear Parabolic {PDEs} with Boundary and In-domain Disturbances}
%
%
%

\author{Jun~Zheng$^{1,2}$\thanks{$^{1}$Department of Basic Course, Southwest Jiaotong University, Emeishan, Sichuan, P. R. of China 614202}
       \thanks{$^{2}$School of Mathematics, Southwest Jiaotong UniversityChengdu, Sichuan, P. R. of China 611756
       {\tt\small zhengjun2014@aliyun.com}}and Guchuan~Zhu$^{3}$,~\IEEEmembership{Senior~Member,~˜IEEE}\thanks{$^{3}$Department of Electrical Engineering, Polytechnique Montr\'{e}al,  P.O. Box 6079, Station Centre-Ville, Montreal, QC, Canada H3T 1J4
       {\tt\small guchuan.zhu@polymtl.ca}}


}

%
%

\markboth{Manuscript Submitted to IEEE Transactions on Automatic Control}%
{Shell \MakeLowercase{\textit{et al.}}: }
%




\maketitle

\begin{abstract}
This paper addresses input-to-state stability (ISS) properties with respect to boundary and in-domain disturbances for a class of semi-linear partial differential equations (PDEs) subject to Dirichlet boundary conditions. The developed approach is a combination of the method of De~Giorgi iteration and the technique of Lyapunov functionals by adequately splitting the original problem into two subsystems. The ISS in $L^\infty$-norm for a 1-$D$ transport equation and the ISS in $L^2$-norm for Burgers' equation have been established using this method. As an application of the main result for the 1-$D$ transport equation, a study on ISS properties in $L^\infty$-norm of a 1-$D$ transport equation with anti-stable term under boundary feedback control is carried out. This is the first time that the method of De~Giorgi iteration is introduced in the ISS theory for infinite dimensional systems, and the developed method can be generalized for tackling some problems on multidimensional spatial domains and be applied to a wider class of nonlinear parabolic PDEs.
\end{abstract}


%
\IEEEpeerreviewmaketitle

\section{Introduction}\label{Sec: Introduction}
Extending the theory of input-to state stability (ISS), which was originally developed for finite-dimensional nonlinear systems \cite{Sontag:1989,Sontag:1990}, to infinite dimensional systems has received a considerable attention in the recent literature. In particular, there are significant progresses on the establishment of ISS estimates with respect to (w.r.t.) disturbances \cite{Argomedo:2013, Argomedo:2012, Dashkovskiy:2010, Dashkovskiy:2013, Dashkovskiy:2013b, jacob2016input, Karafyllis:2014, Karafyllis:2016, Karafyllis:2016a, Karafyllis:2017, Logemann:2013, Mazenc:2011, Prieur:2012, Tanwani:2017,Zheng:2017} for different types of partial differential equations (PDEs).

It is noticed that most of the earlier work on this topic dealt with disturbances distributed over the domain. It was demonstrated that the method of Lyapunov functionals is a well-suited tool for dealing with a wide range of problems of this category. Moreover, it is shown in \cite{Argomedo:2012} that the method of Lyapunov functionals can be readily applied to some systems with boundary disturbances by transforming the later ones to a distributed disturbance. {However, ISS estimates obtained by such a method may include time derivatives of boundary disturbances, which is not strictly in the original form of ISS formulation.

The problems with disturbances acting on the boundaries usually lead to a formulation involving unbounded operators. It is shown in \cite{jacob2016input,Jacob:2017} that for a class of linear PDEs, the exponential stability plus a certain admissibility implies the ISS and iISS (integral input-to-state stability \cite{jacob2016input,Sontag:1998}) w.r.t. boundary disturbances. However, it may be difficult to assess this property for nonlinear PDEs. To resolve this concern while not invoking unbounded operators in the analysis, it is proposed in \cite{Karafyllis:2016,Karafyllis:2016a,karafyllis2017siam} to derive the ISS property directly from the estimates of the solution to the considered PDEs using the method of spectral decomposition and finite-difference. ISS in $L^2$-norm and in weighted $L^\infty$-norm for PDEs with a Sturm-Liouville operator is established by applied this method in \cite{Karafyllis:2016,Karafyllis:2016a,karafyllis2017siam}. However, spectral decomposition and finite-difference schemes  may involve heavy computations for nonlinear PDEs or problems on multidimensional spatial domains. It is introduced in \cite{Mironchenko:2017} a monotonicity-based method for studying the ISS of nonlinear parabolic equations with boundary disturbances. It is shown that with the monotonicity, the ISS of the original nonlinear parabolic PDE with constant boundary disturbances is equivalent to the ISS of a closely related nonlinear parabolic PDE with constant distributed disturbances and zero boundary conditions. As an application of this method, the ISS properties in $L^p$-norm ($p>2$) for some linear parabolic systems have been established.
In a recent work \cite{Zheng:2017}, the classical method of Lyapunov functionals is applied to establish ISS properties in $L^2$-norm w.r.t. boundary disturbances for a class of semi-linear parabolic PDEs. Some technical inequalities have been developed, which allows dealing directly with items on the boundary points in proceeding on ISS estimates. The result of \cite{Zheng:2017} shows that the method of Lyapunov functionals is still effective in obtaining ISS properties for some linear and nonlinear PDEs with Neumann or Robin boundary conditions. However, the technique used may not be suitable for problems with Dirichlet boundary conditions.

The present work deals with a complementary setting compared to that considered in \cite{Zheng:2017} in the sense that the problem is subject to Dirichlet boundary conditions. The method developed in this paper consists first in splitting the original system into two subsystems: a boundary disturbance-free system with homogenous boundary condition and non-zero initial condition, and another one with boundary disturbances and zero-initial condition. Note that the in-domain disturbances can be placed in either of these two subsystems. Then, ISS properties in $H^1$-norm or $L^2$-norm for the first system will be established by the method of Lyapunov functionals, and ISS properties in $L^\infty$-norm for the second system will be deduced by the technique of De~Giorgi iteration. Finally, the ISS properties in $L^\infty$-norm or $L^2$-norm for the original system are obtained by combining the ISS properties of the two subsystems. With this method, we established the ISS in $L^\infty$-norm for a 1-$D$ transport PDE and the ISS in $L^2$-norm for Burgers' equation. In addition, a study on ISS properties in $L^\infty$-norm for a 1-$D$ transport equation with anti-stable term is carried out to illustrate the application of this method for systems with boundary feedback control. Note that although the De~Giorgi iteration is a classic method in regularity analysis of elliptic and parabolic PDEs, it is the first time, to the best of our knowledge, that it is introduced in the investigation of ISS properties for PDEs. Moreover, the technique of De~Giorgi iteration may be applicable for certain nonlinear PDEs and problems on multidimensional spatial domains.

The rest of the paper is organized as follows. Section~\ref{Sec: Preliminaries} introduces briefly the technique of De~Giorgi iteration and presents some preliminary inequalities needed for the subsequent development. Section~\ref{Sec: Main results} presents the considered problems and the main results. Detailed development on the establishment of ISS properties for 1-$D$ transport equations and Burgers' equation is given, respectively, in Section~\ref{Sec: Transport PDE} and Section~\ref{Sec: Burgers Eq}. Finally, some concluding remarks are provided in Section~\ref{Sec: Conclusion}.

\section{Preliminaries}\label{Sec: Preliminaries}

\subsection{{De~Giorgi} iteration}\label{Sec: De Giorgi iteration}

De~Giorgi iteration is an important tool for regularity analysis of elliptic and parabolic PDEs. In his famous work on linear elliptic equations published in 1957 \cite{DeGiorgi:1957}, De~Giorgi established local boundedness and H\"{o}lder continuity for functions satisfying certain integral inequalities, known as the De~Giorgi class of functions, which completed the solution of Hilbert's 19$^{\text{th}}$ problem. The same problem has been resolved independently by Nash in 1958 \cite{Nash1958}. It was shown later by Moser that the result of De~Giorgi and Nash can be obtained using a different formulation \cite{Moser1961}. In the literature, this method is often called the De~Giorgi-Nash-Moser theory.

Let $\Omega\subset \mathbb{R}^N(N\geq 1)$ be an open bounded set and $\gamma $ be a constant. The De~Giorgi class $DG^+(\Omega,\gamma)$ consists of functions $u\in W^{1,2}(\Omega)$ which satisfy, for every ball $B_r(y)\subset \Omega$, every $0<r'<r$, and every $k \in \mathbb{R}$, the following Caccioppoli type inequality:
\begin{align*}
\int_{B_{r'}(y)}|\nabla (u-k)_+|^2\text{d}x\leq \frac{\gamma}{(r-r')^2}\int_{B_{r}(y)}| (u-k)_+|^2\text{d}x,
\end{align*}
where $ (u-k)_+=\max\{u-k,0\}$. The class $DG^-(\Omega,\gamma)$ is defined in a similar way. The main idea of De~Giorgi iteration is to estimate $ |A_k|$, the measure of $\{x\in \Omega;u(x)\geq k\}$, and derive $|A_k|=0$ with some $k$ for functions $u$ in De~Giorgi class by using the iteration formula given below.
\begin{lemma}[\cite{Wu2006}]\label{iteration}
Suppose that $\varphi$ is a non-negative decreasing function on $[k_0,+\infty)$ satisfying
\begin{align*}
\varphi(h)\leq \bigg(\frac{M}{h-k}\bigg)^\alpha\varphi^\beta(k),\ \ \forall h>k\geq k_0,
\end{align*}
where $\alpha>0,\beta>1,M\in \mathbb{R}$ are constants. Then there holds
\begin{align*}
\varphi(k_0+d)=0,
\end{align*}
with $d=2^{\frac{\beta}{\beta-1}}M(\varphi(k_0))^{\frac{\beta-1}{\alpha}}$.
\end{lemma}

The method of De~Giorgi iteration can be generalized to some linear parabolic PDEs and PDEs with a divergence form (see, e.g., \cite{DiBenedetto:2010,Wu2006}). However, this method in its original formulation cannot be applied directly in the establishment of ISS properties for infinite dimensional systems. The main reason is that the obtained boundedness of a solution is not expressed by a $\mathcal{K}\mathcal{L}$-function associated with the initial value and time variable, even for linear parabolic PDEs \cite{DiBenedetto:2010,Wu2006}, as usually expected in ISS analysis. To overcome this difficulty, we developed in this work an approach that amounts first to splitting the original problem into two subsystems and then to applying the De~Giorgi iteration together with the technique of Laypunov functionals to obtain the ISS estimates of the solutions expressed in the standard formulation of the ISS theory.

\subsection{Preliminary inequalities}\label{Sec: preliminary results}
We present below some inequalities needed for the subsequent development.
\begin{lemma}[\textbf{Young's \ inequality}]\label{Lemma 1}
 \begin{align*}
ab\leq \frac{a^2}{2\varepsilon}+\frac{\varepsilon b^2}{2},\ \ \ \ \forall a\geq 0,b\geq 0, \varepsilon>0.
\end{align*}
\end{lemma}
\begin{lemma}[\textbf{Poincar\'{e}'s \ inequality} \cite{Zheng:2017}]\label{Lemma 2}
   Suppose that $u\in C^{1}([a,b];\mathbb{R})$ with $u(c_0)=0$ for some $c_0\in[a,b]$, then
$
\|u\|^2\leq  \frac{(b-a)^2}{2}\|u_x\|^2.
$
\end{lemma}

The following results are some variations of Sobolev embedding inequalities.
\begin{lemma}\label{Lemma 3}
   Suppose that $u\in C^{1}([a,b];\mathbb{R})$, then
   \begin{itemize}
\item[(i)]
$u^2(c)\leq \frac{2}{b-a}\|u\|^2+(b-a)\|u_x\|^2,\ \ \ \ \forall c\in[a,b].
$
\item[(ii)]
$\bigg(\int_{a}^b|u|^p\text{d}x\bigg)^{\frac{1}{p}}\leq (b-a)^{\frac{1}{p}}\bigg(\frac{2}{b-a}\|u\|^2+(b-a)\|u_x\|^2\bigg)^{\frac{1}{2}},\ \ \ \ \forall p\geq 1.
$
\end{itemize}
\end{lemma}
\begin{IEEEproof}
The inequality (i) was proven in \cite{Zheng:2017}. We only prove (ii). Indeed, by (i) we have
\begin{align*}
\bigg(\int_{a}^b|u|^p\text{d}{x}\bigg)^{\frac{1}{p}}&\leq \bigg(\int_{a}^b\max_{x\in[a,b]}|u|^p\text{d}{x}\bigg)^{\frac{1}{p}}\\
&=  (b-a)^{\frac{1}{p}}\max_{x\in[a,b]}|u|\\
&\leq
(b-a)^{\frac{1}{p}}\bigg(\frac{2}{b-a}\|u\|^2+(b-a)\|u_x\|^2\bigg)^{\frac{1}{2}}.
\end{align*}
\end{IEEEproof}
\begin{remark} It is easy to see that Lemma~\ref{Lemma 2} and Lemma~\ref{Lemma 3} hold for $u\in H^1((a,b);\mathbb{R})$ (with $u(c_0)=0$ in Lemma~\ref{Lemma 2}).
\end{remark}
\begin{lemma}[\textbf{Gronwall's \ inequality} \cite{Zheng:2017}]
Suppose that $y:\mathbb{R}_{\geq 0}\rightarrow\mathbb{R}$ is absolutely continuous on $[0,T]$ for any $T>0$ and satisfies for a.e. $t \geq 0$ the following differential inequality
\begin{align*}
\frac{\text{d}y}{\text{d}t}(t)\leq g(t)y(t)+h(t),
\end{align*}
where $g,h\in L^1([0,T];\mathbb{R})$ for any $T > 0$. Then for all $t \in \mathbb{R}_{\geq 0}$,
\begin{align*}
y(t) \leq
y(0)e^{\int_{0}^{t}g(s)\text{d}s}+\int_{0}^{t}h(s) e^{\int_{s}^{t}g(s)\text{d}s}\text{d}s.
\end{align*}
\end{lemma}

\section{Problem Formulation and Main Results}\label{Sec: Main results}

\subsection{Problem formulation and well-posedness analysis}\label{Sec: Problem formulations}
In this paper, we address ISS properties for two classes of PDEs with Dirichlet boundary conditions:

(i) 1-$D$ transport PDE:
\begin{align}\label{++1}
 u_t-\mu u_{xx}+mu_x+nu=f(x,t) \ \ \ \ \text{in}\ \  (0,1)\times \mathbb{R}_{\geq 0},
\end{align}
where $\mu>0$ and $m,n\in \mathbb{R}$ are constants.


(ii) Burgers' equation:
\begin{align}\label{++28}
u_t-\mu u_{xx}+\nu uu_x=f(x,t),\ \ \text{in}\  (0,1)\times \mathbb{R}_{\geq 0},
\end{align}
where $\mu>0$ and $\nu>0$ are constants.

The boundary condition and initial condition for \eqref{++1} and \eqref{++28}  are given by:
\begin{align}
 &u(0,t)=0,u(1,t)=d(t),\label{++2}\\
  &u(x,0)=u_0(x),\label{++3}
\end{align}
where $d(t)$ is the disturbance on the boundary. In general, it can represent actuation and sensing errors. The function $f(x,t)$ in \eqref{++1} and \eqref{++28} can be seen as the disturbance distributed over the domain.

 In order to present the well-posedness of equation \eqref{++1} (or \eqref{++28}) with conditions \eqref{++2} and \eqref{++3}, we consider the following class of semi-linear PDEs:
\begin{subequations}\label{semi-linear PDE}
\begin{align}
&u_t-\mu u_{xx}=g(x,t,u,u_x),\\
 &u(0,t)=0,u(1,t)=d(t),\\
  &u(x,0)=u_0(x),
\end{align}
\end{subequations}

Suppose that $g:(0,1)\times\mathbb{R}_{\geq 0}\times \mathbb{R}\times\mathbb{R}\rightarrow \mathbb{R}$ is locally Lipschitz continuous in all its variables and satisfies
\begin{subequations}\label{AF}
\begin{align}
|g(x,t,u,p)|&\leq \rho(t,|u|)(1+|p|^\gamma),\label{AF-1}\\
|g(x,t,u,p)&-g(t,x,u,q)| \nonumber \\
           &\leq \rho(t,|u|)(1+|p|^{\gamma-1}+|q|^{\gamma-1})(|p-q|),\label{AF-2}\\
|g(x,t,u,p)&-g(t,x,v,p)| \nonumber \\
           &\leq \rho(t,|u|+|v|)(1+|p|^\gamma)(|u-v|),\label{AF-3}
\end{align}
\end{subequations}
 for almost all $x\in (0,1)$ and all $ t\in\mathbb{R}_{\geq 0} , u\in\mathbb{R},v\in\mathbb{R},p\in\mathbb{R},q\in\mathbb{R}$, where $\rho:\mathbb{R}_{\geq 0}\times\mathbb{R}\rightarrow\mathbb{R}_{\geq 0}$ is continuous in both variables and is monotonously increasing in the second argument, and $\gamma\in [1,3) $ is a constant.


\begin{proposition}\label{Proposition 1}
Suppose that $g$ satisfies the structural conditions \eqref{AF}, $f(x,t)\in C^1((0,1)\times \mathbb{R}_{\geq 0})$ and $d\in C^2(\mathbb{R}_{\geq 0})$. If
$u_0-d(0)\psi\in H^2(0,1)\cap H^1_0(0,1)$ for some $\psi\in C^3([0,1])$ and $d(0)$ with $d(0)=u_0(1)$, then
there exists a unique solution $u\in C((0,1);H^2(0,1))\cap C^{1}((0,1);L^2(0,1))$ of \eqref{semi-linear PDE}.
\end{proposition}
\begin{IEEEproof}
The proof is based on the technique of lifting (see, e.g., \cite[Theorem~4]{Zheng:2017}) and the theory of Lipschitz perturbations for linear evolution equations (see, e.g., \cite[\S6.3, Chap.~6]{Pazy:1983}). As the proof can be proceeded with a standard procedure, it is omitted and we refer to \cite{Zheng:2017} for details (see also \cite[Proposition 7, \S4.4]{Kuttler:2011}).
\end{IEEEproof}

Throughout this paper, we always assume that $f(x,t)\in C^1([0,1]\times \mathbb{R}_{\geq 0})$, $d\in C^2(\mathbb{R}_{\geq 0})$, $u_0-d(0)\psi\in H^2(0,1)\cap H^1_0(0,1)$ for some $\psi\in C^3([0,1])$ and $d(0)$ with $d(0)=u_0(1)$, and $\mu>0$ is a constant. Unless specially stated, we always take $u\in C((0,1);H^2(0,1))\cap C^{1}((0,1);L^2(0,1))$ as the solution of \eqref{++1} or \eqref{++28}. For notational simplicity, we always denote $\|\cdot\|_{L^{2}(0,1)}$ by $\|\cdot\|$.

\subsection{Main results}
Let $\mathcal {K}=\{\gamma : \mathbb{R}_{\geq 0} \rightarrow \mathbb{R}_{\geq 0}|\ \gamma(0)=0,\gamma$ is continuous, strictly increasing$\}$; $\mathcal {K}_{\infty}=\{\theta \in \mathcal {K}|\ \lim\limits_{s\rightarrow\infty}\theta(s)=\infty\}$; $\mathcal {L}=\{\gamma : \mathbb{R}_{\geq 0}\rightarrow \mathbb{R}_{\geq 0}|\ \gamma$ is continuous, strictly decreasing, $\lim\limits_{s\rightarrow\infty}\gamma(s)=0\}$; $\mathcal {K}\mathcal {L}=\{\beta : \mathbb{R}_{\geq 0}\times \mathbb{R}_{\geq 0}\rightarrow \mathbb{R}_{\geq 0}|\ \beta(\cdot,t)\in \mathcal {K}, \forall t \in \mathbb{R}_{\geq 0}$, and $\beta(s,\cdot)\in \mathcal {L}, \forall s \in \mathbb{R}_{\geq 0}\}$.

Let $u_0$ be the initial state of \eqref{++1} (or \eqref{++28}) with ~\eqref{++2} and \eqref{++3} in a certain space $\mathcal{H}$ with norm $\|\cdot\|_\mathcal{H}$.

\begin{definition}
System~\eqref{++1} (or \eqref{++28}) is said to be input-to-state stable (ISS) in $L^\infty-$norm (or $L^2-$norm) w.r.t. the boundary disturbance $d(t)$ and the in-domain disturbance $f(x,t)$, if there exist functions $\beta\in \mathcal {K}\mathcal {L}$ and $ \gamma_1, \gamma_2,\in \mathcal {K}$ such that the solution of \eqref{++1} (or \eqref{++28}) satisfies
\begin{align}\label{Eq: ISS def}
\begin{split}
     \max_{x\in[0,1]}|u(x,t)|\leq & \beta( \|{u_0}\|_\mathcal{H},t)
      +\gamma_1\left(\max_{s\in [0,t]}|d(s)|\right) \\
      &+\gamma_2\left(\max_{(x,s)\in [0,1]\times [0,t]}|f(x,s)|\right),\ \forall t\geq 0,
\end{split}
\end{align}
or
\begin{align}\label{Eq: ISS def2}
\begin{split}
\|u\|\leq & \beta( \|{u_0}\|_\mathcal{H},t)+\gamma_1\left(\max_{s\in [0,t]}|d(s)|\right)\\
          &+\gamma_2\left(\max_{(x,s)\in [0,1]\times [0,t]}|f(x,s)|\right),\ \forall t\geq 0.
\end{split}
\end{align}
Moreover, System~\eqref{++1} (or \eqref{++28}) is said to be exponential input-to-state stable (EISS) w.r.t. the boundary disturbance $d(t)$ and the in-domain disturbance $f(x,t)$, if there exist $\beta'\in \mathcal {K}_{\infty}$ and a constant $\lambda > 0$ such that  $\beta( \|{u_0}\|_\mathcal{H},t) \leq \beta'(\|{u_0}\|_\mathcal{H})e^{-\lambda t}$ in \eqref{Eq: ISS def} (or \eqref{Eq: ISS def2}).
\end{definition}
\begin{definition}
System~\eqref{++1} (or \eqref{++28}) is said to be ISS in $L^\infty$-norm (or $L^2$-norm) w.r.t. the boundary disturbance $d(t)$ and iISS w.r.t. the in-domain disturbance $f(x,t)$, if there exist functions $\beta\in \mathcal {K}\mathcal {L},\theta\in \mathcal {K}_{\infty} $, and $\gamma_1 ,\gamma_2 \in \mathcal {K}$ such that the solution of \eqref{++1} (or \eqref{++28}) satisfies
\begin{align}\label{Eq: iISS def}
\begin{split}
     \max_{x\in[0,1]}|u(x,t)|\leq & \beta( \|{u_0}\|_\mathcal{H},t)
      +\gamma_1\left(\max_{s\in [0,t]}|d(s)|\right) \\
      &+\theta\bigg(\!\!\int_{0}^t\!\!\gamma_2(\|f(\cdot,s)\|)\text{d}s\bigg),\ \forall t\geq 0,
\end{split}
\end{align}
or
\begin{align}\label{Eq: iISS def2}
\begin{split}
\|u\| \leq & \beta( \|{u_0}\|_\mathcal{H},t)
      +\gamma_1\left(\max_{s\in [0,t]}|d(s)|\right) \\
      &+\theta\bigg(\!\!\int_{0}^t\!\!\gamma_2(\|f(\cdot,s)\|)\text{d}s\bigg),\ \forall t\geq 0.
\end{split}
\end{align}

Moreover, System~\eqref{++1} (or \eqref{++28}) is said to be EISS w.r.t. the boundary disturbance $d(t)$ and exponential integral input-to-state stable (EiISS) w.r.t. the in-domain disturbance $f(x,t)$, if there exist $\beta'\in \mathcal {K}_{\infty}$ and a constant $\lambda > 0$ such that  $\beta( \|{u_0}\|_\mathcal{H},t) \leq \beta'(\|{u_0}\|_\mathcal{H})e^{-\lambda t}$ in \eqref{Eq: iISS def} (or \eqref{Eq: iISS def2}).
\end{definition}

The ISS properties w.r.t. boundary and in-domain disturbances for systems \eqref{++1} and \eqref{++28} are stated in the following theorems.

\begin{theorem}\label{Theorem 6}
Assume $m,n\in \mathbb{R}$ satisfying $\frac{m^2}{4\mu}+n\geq 0$. Then the following statements hold true.
\begin{itemize}
\item[(i)]
System \eqref{++1} with \eqref{++2} and \eqref{++3} is EISS w.r.t. the boundary and in-domain disturbances with the following estimate for any $t>0$:
\begin{align*}
&\max_{x\in[0,1]}|u(x,t)| \notag\\
&\;\;\; \leq e^{\frac{|m|}{\mu}}\sqrt{2+\frac{m^2}{2\mu^2}} \|u_0\|_{H^1(0,1)}e^{-(\frac{m^2}{4\mu}+n+2\mu) t} \notag\\
&\;\;\;\;\;\;  +e^{\frac{|m|}{\mu}}\max\limits_{s\in [0,t]}|d(s)|\notag\\
&\;\;\;\;\;\;  +\frac{1}{\mu}e^{\frac{|m|}{\mu}}2^{\frac{5p-8}{2p-4}}\max\limits_{(x,s)\in[0,1]\times [0,t]}|f(x,s)|,\ \ \forall p>2.
\end{align*}
\end{itemize}
\begin{itemize}
\item[(ii)]
System \eqref{++1} with \eqref{++2} and \eqref{++3} is EISS w.r.t. the boundary disturbance and EiISS w.r.t. the in-domain disturbance with the following estimate:
\begin{align*}
&\max_{x\in[0,1]}|u(x,t)|\notag\\
&\;\;\; \leq e^{\frac{|m|}{\mu}}\sqrt{2+\frac{m^2}{2\mu^2}}\|u_0\|_{H^1(0,1)}e^{-(\frac{m^2}{4\mu}+n+2\mu-\frac{\varepsilon}{2}) t} \notag\\
&\;\;\;\;\;\; +e^{\frac{|m|}{\mu}}\max\limits_{s\in [0,t]}|d(s)|\notag\\
&\;\;\;\;\;\; +e^{\frac{|m|}{\mu}}\sqrt{\frac{3}{\varepsilon} \int_{0}^t\|f(\cdot,s)\|^2\text{d}s},
\end{align*}
for any $t>0$ and all $\varepsilon \in \left(0, \displaystyle\frac{m^2}{2\mu}+2n+4\mu\right)$.
\end{itemize}
\end{theorem}



\begin{theorem} \label{Theorem 11}
The following statements hold true.
\begin{itemize}
\item[(i)] Suppose that $\nu>0$ and $\max\limits_{t\in \mathbb{R}_{\geq 0}} |d(t)|<\frac{\mu}{\nu}$. System \eqref{++28} with \eqref{++2} and \eqref{++3} is EISS w.r.t. the boundary disturbance and EiISS w.r.t. the in-domain disturbance with the following estimate for any $t>0$:
 \begin{align*}
\|u(\cdot,t)\|^2 \leq& 2\|u_0\|^2 e^{-(\mu -\varepsilon)t}+2\max\limits_{s\in [0,t]}|d(s)|^2 \nonumber \\
                     &+\frac{2}{\varepsilon} \int_{0}^t\|f(\cdot,s)\|^2\text{d}s,\ \forall\varepsilon\in (0,\mu).
\end{align*}
\end{itemize}
\begin{itemize}
\item[(ii)]Suppose that $\nu>0$ and for some $p \in (2,+\infty)$ $\max\limits_{t\in \mathbb{R}_{\geq 0}} |d(t)| + \frac{1}{\mu}2^{\frac{5p-8}{2p-4}}\max\limits_{(x,s)\in[0,1]\times \mathbb{R}_{\geq 0}}|f(x,s)| < \frac{\mu}{\nu}$. System \eqref{++28} with \eqref{++2} and \eqref{++3} is EISS w.r.t. the boundary and in-domain disturbances with the following estimate for any $t>0$:
 \begin{align*}
\|u(\cdot,t)\|^2\leq & 2e^{-\mu t}\|u_0\|^2 +4\max\limits_{s\in [0,t)} |d(s)|^2 \nonumber \\
                     &+\frac{1}{\mu^2}2^{2+\frac{5p-8}{2p-4}}\max\limits_{(x,s)\in[0,1]\times [0,t]}|f(x,s)|^2.
\end{align*}
\end{itemize}
\end{theorem}

\section{ISS Properties for 1-$D$ Transport PDE}\label{Sec: Transport PDE}
\subsection{ISS properties for the split subsystems}
In order to establish the ISS properties in $L^\infty$-norm for the 1-$D$ transport PDE \eqref{++1} with boundary and initial conditions \eqref{++2} and \eqref{++3}, we consider the following two systems and establish their stabilities:
\begin{subequations}\label{++7}
 \begin{align}\label{++7.1}
 &w_t-\mu w_{xx}+mw_x+nw=0,\\
 &w(0,t)=w(1,t)=0,\label{++7.2}\\
  &w(x,0)=u_0(x),\label{++7.3}
\end{align}
\end{subequations}
and
\begin{subequations}\label{++8}
 \begin{align}\label{++8.1}
 &v_t-\mu v_{xx}+mv_x+nv=f(x,t),\\
 &v(0,t)=0,v(1,t)=d(t),\label{++8.2}\\
  &v(x,0)=0.\label{++8.3}
\end{align}
\end{subequations}
The existence and uniqueness of a solution $w$ of \eqref{++7} and $v$ of \eqref{++8} is guaranteed by Proposition \ref{Proposition 1}. Moreover, the solution of \eqref{++1} with boundary and initial conditions \eqref{++2} and \eqref{++3} is given by $u=w+v$.

For \eqref{++7}, we have the following results.
\begin{theorem}\label{Theorem 7}
Let $w\in  C((0,1);H^2(0,1))\cap C^{1}((0,1);L^2(0,1))$ be the unique solution of \eqref{++7}.  For every $t>0$, there holds
\begin{align*}
\|w\|_{L^2(0,1)}^2 \leq  &e^{\frac{|m|}{\mu}}\|u_0\|_{L^2(0,1)}^2e^{-2(\frac{m^2}{4\mu}+n+2\mu)t},\notag\\
\|w_x\|_{L^2(0,1)}^2\leq &e^{\frac{2|m|}{\mu}} \bigg(  \frac{3m^2}{2\mu^2}\|u_0\|_{L^2(0,1)}^2+ 4\| u_{0x}\|_{L^2(0,1)}^2 \bigg) \nonumber \\
                         &\times e^{-2(\frac{m^2}{4\mu}+n+2\mu)t}.\notag
\end{align*}
\end{theorem}
\begin{theorem}\label{Corollary 8}
Let $w\in  C((0,1);H^2(0,1))\cap C^{1}((0,1);L^2(0,1))$ be the unique solution of \eqref{++7}. For every $t>0$, there holds
\begin{align*}
\max_{x\in[0,1]}|w(x,t)|\leq e^{\frac{|m|}{\mu}}\sqrt{2+\frac{m^2}{2\mu^2}}\|u_0\|_{H^1(0,1)} e^{-(\frac{m^2}{4\mu}+n+2\mu)t}.\notag
\end{align*}
\end{theorem}
For \eqref{++8}, we have the following result.
\begin{theorem}\label{Theorem 9}
Let $v\in  C((0,1);H^2(0,1))\cap C^{1}((0,1);L^2(0,1))$ be the unique solution of \eqref{++8}. For every $t>0$ and all $p>2$, there holds:
\begin{align*}
&\max\limits_{(x,s)\in[0,1]\times [0,t]}|v(x,s)| \nonumber \\
&\;\;\;\leq e^{\frac{|m|}{\mu}}\bigg(\max\limits_{s\in [0,t]}|d(s)|+\frac{ 1 }{\mu}2^{\frac{5p-8}{2p-4}} \max\limits_{(x,s)\in[0,1]\times [0,t]}|f(x,s)|\bigg).
\end{align*}
\end{theorem}

\begin{IEEEproof}[Proof of Theorem~\ref{Theorem 7}]
For the convenience of computations, let $w(x,t)=\widetilde{w}(x,t)e^{-\frac{mx}{2\mu}}$, we can transform \eqref{++7} to the following PDE:
\begin{subequations}\label{++9}
 \begin{align}\label{++9.1}
 &\widetilde{w}_t-\mu \widetilde{w}_{xx}+\widetilde{n}\widetilde{w}=0,\\
 &\widetilde{w}(0,t)=\widetilde{w}(1,t)=0,\label{++9.2}\\
  &\widetilde{w}(x,0)=\widetilde{w}_0,\label{++9.3}
\end{align}
\end{subequations}
where $\widetilde{n}=\frac{m^2}{4\mu}+n$ and $\widetilde{w}_0=e^{\frac{mx}{2\mu}}u_0(x) $.

Multiplying \eqref{++9} by $\widetilde{w}$ and integrating over $(0,1)$, one may get
\begin{align*}
 \int_{0}^{1}\widetilde{w}_t\widetilde{w}\text{d}x+\mu\int_{0}^{1}\widetilde{w}_x^2\text{d}x+\widetilde{n}\int_{0}^{1}\widetilde{w}^2\text{d}x=0,
\end{align*}
which is
\begin{align*}
 \frac{1}{2}\frac{d}{dt}\|\widetilde{w}\|^2+\mu\|\widetilde{w}_x\|^2+\widetilde{n}\|\widetilde{w}\|^2=0.
\end{align*}
Note that $ \|\widetilde{w}_x\|^2\geq 2\|\widetilde{w}\|^2$. It follows
\begin{align*}
 \frac{d}{dt}\|\widetilde{w}\|^2+2(2\mu+\widetilde{n})\|\widetilde{w}\|^2\leq 0.
\end{align*}
By Gronwall's inequality, we may obtain
\begin{align}
 \|\widetilde{w}(\cdot,t)\|^2&\leq \|\widetilde{w}_0\|^2e^{-2(2\mu+\widetilde{n})t}\notag\\
 &\leq \|e^{\frac{mx}{2\mu}}u_0\|^2e^{-2(2\mu+\widetilde{n})t}\notag\\
 &\leq e^{\frac{|m|}{\mu}}\|u_0\|^2e^{-2(2\mu+\widetilde{n})t},\label{++10}
\end{align}
which yields
\begin{align*}
\|w(\cdot,t)\|^2 =&\|e^{-\frac{mx}{2\mu}}\widetilde{w}(\cdot,t)\|^2\leq \|\widetilde{w}(\cdot,t)\|^2 \\
                 \leq &  e^{\frac{|m|}{\mu}}\|u_0\|^2e^{-2(2\mu+\widetilde{n})t}.
\end{align*}

Multiplying \eqref{++9} by $\widetilde{w}_{xx}$ and integrating over $(0,1)$, one may get
\begin{align*}
 \int_{0}^{1}\widetilde{w}_{tx}\widetilde{w}_x\text{d}x+\mu\int_{0}^{1}\widetilde{w}_{xx}^2\text{d}x+\widetilde{n}\int_{0}^{1}\widetilde{w}_x^2\text{d}x=0,
\end{align*}
which is
\begin{align*}
 \frac{1}{2}\frac{d}{dt}\|\widetilde{w}_x\|^2+\mu\|\widetilde{w}_{xx}\|^2+\widetilde{n}\|\widetilde{w}_x\|^2=0.
\end{align*}
Note that $\widetilde{w}(0,t)=\widetilde{w}(1,t)=0 $ and $\widetilde{w}$ is $C^1$-continuous in $x$, by Rolle's theorem \cite{Stromberg1981}, there exits $c_0\in(0,1)$ such that $\widetilde{w}_x(c_0,t)=0$. Thus, Lemma~\ref{Lemma 2} gives $\|\widetilde{w}_{xx}\|^2\geq 2\|\widetilde{w}_{x}\|^2$.
It follows
\begin{align*}
 \frac{d}{dt}\|\widetilde{w}_x\|^2+2(2\mu+\widetilde{n})\|\widetilde{w}_x\|^2\leq 0.
\end{align*}
By Gronwall's inequality, we may obtain
\begin{align}
 \|\widetilde{w}_x(\cdot,t)\|^2&\leq \|\widetilde{w}_{0x}\|^2e^{-2(2\mu+\widetilde{n})t}\notag\\
 &\leq \|\frac{m}{2\mu} e^{\frac{mx}{2\mu}} u_0+e^{\frac{mx}{2\mu}} u_{0x}\|^2e^{-2(2\mu+\widetilde{n})t}\notag\\
 &\leq e^{\frac{|m|}{\mu}}\bigg(\frac{m^2}{2\mu^2}  \| u_0\|^2+2\| u_{0x}\|^2\bigg)e^{-2(2\mu+\widetilde{n})t} ,\label{++11}
\end{align}
from which and \eqref{++10} it yields
\begin{align*}
\|w_x(\cdot,t)\|^2 =&\|-\frac{m}{2\mu}e^{-\frac{mx}{2\mu}}\widetilde{w}(\cdot,t)+e^{-\frac{mx}{2\mu}}\widetilde{w}_x(\cdot,t)\|^2\notag\\
\leq &  e^{\frac{|m|}{\mu}}\bigg(    \frac{m^2}{2\mu^2} \|\widetilde{w}(\cdot,t)\|^2+2\|\widetilde{w}_x(\cdot,t)\|^2\bigg) \notag\\
\leq & e^{\frac{|m|}{\mu}}\bigg(  e^{\frac{|m|}{\mu}}\frac{m^2}{2\mu^2}\|u_0\|^2e^{-2(2\mu+\widetilde{n})t} \notag\\
     &+2e^{\frac{|m|}{\mu}}\bigg(\frac{m^2}{2\mu^2}  \| u_0\|^2+2\| u_{0x}\|^2\bigg)e^{-2(2\mu+\widetilde{n})t}\bigg) \notag\\
\leq &e^{\frac{2|m|}{\mu}}\bigg( \frac{3m^2}{2\mu^2}\|u_0\|^2+ 4\| u_{0x}\|^2   \bigg)e^{-2(2\mu+\widetilde{n})t}.
\end{align*}
\end{IEEEproof}

\begin{IEEEproof}[Proof of Theorem~\ref{Corollary 8}] Indeed, for every $t>0$, by the definition of $\widetilde{w} $, (i) of Lemma~\ref{Lemma 3}, \eqref{++10} and \eqref{++11}, we have
\begin{align*}
\max_{x\in[0,1]}|w(x,t)|
=& \max_{x\in[0,1]}\big|\widetilde{w}(x,t)e^{-\frac{mx}{2\mu}}\big| \notag\\
\leq &e^{\frac{|m|}{2\mu}}\max\limits_{x\in [0,1]} |\widetilde{w}(x,t)|\notag\\
\leq &e^{\frac{|m|}{2\mu}}\sqrt{2 \|\widetilde{w}(\cdot,t)\|^2+\|\widetilde{w}_x(\cdot,t)\|^2}\notag\\
\leq &e^{\frac{|m|}{2\mu}}\sqrt{e^{\frac{|m|}{\mu}}\bigg(\bigg(2+\frac{m^2}{2\mu^2}\bigg)\|u_0\|^2+  2\| u_{0x}\|^2 \bigg)} \notag\\
     &\times e^{-(\frac{m^2}{4\mu}+n+2\mu)t}\notag\\
\leq &e^{\frac{|m|}{\mu}} \sqrt{2+\frac{m^2}{2\mu^2} }\|u_0\|_{H^1(0,1)} e^{-(\frac{m^2}{4\mu}+n+2\mu)t}.
\end{align*}
\end{IEEEproof}


\begin{IEEEproof}[Proof of Theorem~\ref{Theorem 9}]
Arguing as the proof of Theorem~\ref{Theorem 7}, let $v(x,t)=\widetilde{v}(x,t)e^{-\frac{mx}{2\mu}}$, we can transform \eqref{++8} to the following PDE:
\begin{subequations}\label{++12}
 \begin{align}\label{++12.1}
 &\widetilde{v}_t-\mu \widetilde{v}_{xx}+\widetilde{n}\widetilde{v}=\widetilde{f}(x,t),\\
 &\widetilde{v}(0,t)=0,\widetilde{v}(1,t)=\widetilde{d}(t),\label{++12.2}\\
  &\widetilde{v}(x,0)=0,\label{++12.3}
\end{align}
\end{subequations}
where $\widetilde{f}(x,t)=e^{\frac{mx}{2\mu}}f(x,t)$, $\widetilde{d}(t)=e^{\frac{m}{2\mu}}d(t) $ and $\widetilde{n}=\frac{m^2}{4\mu}+n$.

We show now that for any $t>0$ and all $p>2$, there holds
\begin{align}\label{++13}
&\max\limits_{(x,s)\in[0,1]\times [0,t]} |\widetilde{v}(x,s)| \notag\\
&\;\;\; \leq \max\limits_{s\in [0,t]}|\widetilde{d}(s)| + \frac{1}{\mu}2^{\frac{5p-8}{2p-4}} \max\limits_{(x,s)\in[0,1]\times [0,t]}|\widetilde{f}(x,s)|.
\end{align}
For this aim, we resort to the De~Giorgi iteration, which can be conducted with a standard process (see, e.g., \cite{Wu2006}). However, we provide below the details for the completeness of the development. For any $t>0$, let $Q_t= (0,1)\times (0,t)$ and $W^{1,0}_2(Q_t)=\{\psi\in L^{2}(Q_t); \psi_x\in L^{2}(Q_t)\}$. Let $\mathring{C}(\overline{Q}_t)=\{\psi\in C^\infty( \overline{Q}_t); \psi(0,s)=\psi(1,s)=0, \  \forall s \in (0,t)\}$ and let $\mathring{W}^{1,0}_2(Q_t)$ be the closure of $\mathring{C}(\overline{Q}_t)$ in $ W^{1,0}_2(Q_t)$.

Now for any fixed $t>0$, let $k_0=\max\Big\{\max\limits_{s\in[0,t]}\widetilde{d}(s),0\Big\}$. For any $k\geq k_0$, let $ \eta(x,s)=(\widetilde{v}(x,s)-k)_+\chi_{[t_1,t_2]}(
s)$, which belongs to $\mathring{W}^{1,0}_2(Q_t)$, where $\chi_{[t_1,t_2]}(t) $ is the character function on $[t_1,t_2]$ and $0\leq t_1<t_2\leq t$. Multiplying \eqref{++12.1} by $\eta$, we get
\begin{align}\label{+16}
&\int_{0}^t\int_{0}^1(\widetilde{v}-k)_t(\widetilde{v}-k)_+\chi_{[t_1,t_2]}(s)\text{d}x\text{d}s \notag\\
&+\mu\int_{0}^t\int_{0}^1\chi_{[t_1,t_2]}(s)|((\widetilde{v}-k)_+)_x |^2\text{d}x\text{d}s\notag\\
&+\widetilde{n}\int_{0}^t\int_{0}^1\widetilde{v}(\widetilde{v}-k)_+\chi_{[t_1,t_2]}(s) \text{d}x\text{d}s \notag\\
=&\int_{0}^t\int_{0}^1\widetilde{f}(\widetilde{v}-k)_+\chi_{[t_1,t_2]}(s) \text{d}x\text{d}s.
\end{align}
Let $I_k(s)=\int_{0}^1(\widetilde{v}(x,s)-k)_+^2\text{d}x$, which is absolutely continuous on $[0,t]$. Suppose that $I_k(t_0)=\max\limits_{s\in[0,t]}I_k(s)$ with some $t_0\in [0,t]$. Due to $I_k(0)=0$ and $I_k(s)\geq 0 $, one may assume that $t_0>0$ without loss of generality.

Note that
$
\widetilde{n}\widetilde{v}(\widetilde{v}-k)_+\chi_{[t_1,t_2]}(s)\geq 0.
$
For $ \varepsilon>0$ small enough, choosing $t_1=t_0-\varepsilon$ and $t_2=t_0 $, it follows
\begin{align*}
&\frac{1}{2\varepsilon}\int_{t_0-\varepsilon}^{t_0}\frac{d}{dt}\int_{0}^1(\widetilde{v}-k)_+^2\text{d}x\text{d}s \notag\\
&+\frac{\mu}{\varepsilon}\int_{t_0-\varepsilon}^{t_0}\int_{0}^1|((\widetilde{v}-k)_+)_x |^2\text{d}x\text{d}s \notag\\
\leq & \frac{1}{\varepsilon}\int_{t_0-\varepsilon}^{t_0}\int_{0}^1|\widetilde{f}|(\widetilde{v}-k)_+\text{d}x\text{d}s.
\end{align*}
Note that
 \begin{align*}
\frac{1}{2\varepsilon}\int_{t_0-\varepsilon}^{t_0}\frac{d}{dt}\int_{0}^1(\widetilde{v}-k)_+^2\text{d}x\text{d}s=\frac{1}{2\varepsilon}( I_k(t_0)-I_k(t_0-\varepsilon))\geq 0.
\end{align*}
We have
\begin{align*}
     & \frac{\mu}{\varepsilon}\int_{t_0-\varepsilon}^{t_0}\int_{0}^1|((\widetilde{v}-k)_+)_x |^2\text{d}x\text{d}s \notag\\
\leq & \frac{1}{\varepsilon}\int_{t_0-\varepsilon}^{t_0}\int_{0}^1|\widetilde{f}|(\widetilde{v}-k)_+\text{d}x\text{d}s.
\end{align*}
Letting $ \varepsilon\rightarrow 0^+$, we get
\begin{align*}
&\mu\int_{0}^1|((\widetilde{v}(x,t_0)-k)_+)_x |^2\text{d}x \notag\\
\leq & \int_{0}^1|\widetilde{f}(x,t_0)|(\widetilde{v}(x,t_0)-k)_+\text{d}x.
\end{align*}

Due to $(\widetilde{v}(0,t_0)-k)_+ =(\widetilde{v}(1,t_0)-k)_+=0$, we deduce by (ii) of Lemma~\ref{Lemma 3} and Lemma~\ref{Lemma 2} that for all $p>2$,
\begin{align*}
     &\bigg(\int_{0}^1|((\widetilde{v}(x,t_0)-k)_+) |^p\text{d}x\bigg)^{\frac{2}{p}} \notag\\
\leq & 2 \int_{0}^1|((\widetilde{v}(x,t_0)-k)_+)_x |^2\text{d}x \notag\\
\leq & \frac{2}{\mu} \int_{0}^1|\widetilde{f}(x,t_0)|(\widetilde{v}(x,t_0)-k)_+\text{d}x.
\end{align*}

Let $A_{k}(s)=\{x\in (0,1);\widetilde{v}(x,s)>k\}$ and $\varphi_{k}=\max\limits_{s\in[0,t]}|A_{k}(s)|$, where $|B|$ denotes the 1-dimensional Lebesgue measure of a set $B\subset(0,1)$. Then we have
\begin{align*}
&\bigg(\int_{A_{k}(t_0)}|((\widetilde{v}(x,t_0)-k)_+) |^p\text{d}x\bigg)^{\frac{2}{p}} \notag\\
\leq & \frac{2}{\mu}   \int_{A_{k}(t_0)}|\widetilde{f}(x,t_0)|(\widetilde{v}(x,t_0)-k)_+\text{d}x.
\end{align*}
By H\"{o}lder inequality, it follows
\begin{align*}
&\bigg(\int_{A_{k}(t_0)}|((\widetilde{v}(x,t_0)-k)_+) |^p\text{d}x\bigg)^{\frac{2}{p}} \notag\\
\leq &\frac{2}{\mu} \bigg(\int_{A_{k}(t_0)}|(\widetilde{v}(x,t_0)-k)_+^p\text{d}x\bigg)^{\frac{1}{p}}\bigg(\int_{0}^1|\widetilde{f}(x,t_0)|^q\text{d}x\bigg)^{\frac{1}{q}},
\end{align*}
where $\frac{1}{p}+\frac{1}{q}=1$.
Thus
\begin{align}\label{++14}
&\bigg(\int_{A_{k}(t_0)}|((\widetilde{v}(x,t_0)-k)_+) |^p\text{d}x\bigg)^{\frac{1}{p}}\notag\\
\leq & \frac{2}{\mu} \bigg(\int_{A_{k}(t_0)}|\widetilde{f}(x,t_0)|^q\text{d}x\bigg)^{\frac{1}{q}}\notag\\
\leq & \frac{2}{\mu}  \max_{(x,s)\in [0,1]\times[0,t]}|\widetilde{f}(x,s)||{A_{k}(t_0)}|^{\frac{1}{q}}\notag\\
\leq & \frac{2}{\mu}  \max_{(x,s)\in [0,1]\times[0,t]}|\widetilde{f}(x,s)|\varphi_{k}^{\frac{1}{q}}.
\end{align}
Now for $I_k(t_0)$, we get by H\"{o}lder inequality and \eqref{++14}
\begin{align*}
I_k(t_0)&\leq \bigg(\int_{A_{k}(t_0)}|((\widetilde{v}(x,t_0)-k)_+) |^p\text{d}x \bigg)^{\frac{2}{p}}|{A_{k}(t_0)}|^{\frac{p-2}{p}}\notag\\
&\leq \bigg(\frac{2}{\mu} \max_{(x,s)\in [0,1]\times[0,t]}|\widetilde{f}(x,s)|\bigg)^2\varphi_{k}^{3-\frac{4}{p}}.
\end{align*}
Recalling the definition of $I_k(t_0)$, for any $s\in [0,t]$ we conclude that
\begin{align}
I_k(s)\leq I_k(t_0)\leq \bigg(\frac{2}{\mu} \max_{(x,s)\in [0,1]\times[0,t]}|\widetilde{f}(x,s)|\bigg)^2\varphi_{k}^{3-\frac{4}{p}}.\label{++15}
\end{align}
Note that for any $h>k$ and $s\in [0,t]$ there holds
\begin{align*}
I_k(s)\geq \int_{A_{h}(t_0)}(\widetilde{v}(x,t_0)-k)_+^2 \text{d}x\geq (h-k)^2|A_h(s)|.
\end{align*}
Then we infer from \eqref{++15} that
\begin{align*}
(h-k)^2\varphi_h\leq \bigg(\frac{2}{\mu} \max_{(x,s)\in [0,1]\times[0,t]}|\widetilde{f}(x,s)|\bigg)^2\varphi_{k}^{3-\frac{4}{p}},
\end{align*}
which is
\begin{align*}
\varphi_h\leq \left(\frac{2}{\mu}\frac{\max\limits_{(x,s)\in [0,1]\times[0,t]}|\widetilde{f}(x,s)|}{h-k}\right)^2\varphi_{k}^{3-\frac{4}{p}}.
\end{align*}
As $p>2$, we have $ 3-\frac{4}{p}>1$. By Lemma~\ref{iteration}, we obtain
\begin{align*}
\varphi_{k_0+d}=\max_{s\in[0,t]}|A_{k_0+d}|=0,
\end{align*}
where $d=2^{\frac{3p-4}{2p-4}}\frac{2}{\mu} \max\limits_{(x,s)\in [0,1]\times[0,t]}|\widetilde{f}(x,s)|\varphi_{k_0}^{1-\frac{2}{p}}\leq \frac{1}{\mu}2^{\frac{5p-8}{2p-4}}\max\limits_{(x,s)\in [0,1]\times[0,t]}|\widetilde{f}(x,s)|$.

By the definition of $A_k$, for almost all $(x,s)\in [0,1]\times [0,t]$ there holds
\begin{align*}
\widetilde{v}(x,s)
\leq & k_0+\frac{1}{\mu}2^{\frac{5p-8}{2p-4}}\max\limits_{(x,s)\in [0,1]\times[0,t]}|\widetilde{f}(x,s)|\notag\\
 = &\max\Big\{\max\limits_{s\in[0,t]}\widetilde{d}(s),0\Big\} \notag\\
   &+\frac{1}{\mu}2^{\frac{5p-8}{2p-4}}\max\limits_{(x,s)\in [0,1]\times[0,t]}|\widetilde{f}(x,s)|.
\end{align*}
By continuity of $\widetilde{v}(x,s) $, for every $(x,s)\in [0,1]\times [0,t]$ there holds
\begin{align}
\widetilde{v}(x,s)
\leq & \max\{\max\limits_{s\in[0,t]}\widetilde{d}(s),0\}\notag\\
     &+\frac{1}{\mu}2^{\frac{5p-8}{2p-4}}\max\limits_{(x,s)\in [0,1]\times[0,t]}|\widetilde{f}(x,s)|.\label{++16}
\end{align}
To conclude on the inequality \eqref{++13}, we should prove the lower boundedness of $\widetilde{v}(x,t)$. Indeed, setting $\overline{v}=-\widetilde{v}$, we get
\begin{align*}
 &\overline{v}_t-\mu \overline{v}_{xx}+\widetilde{n}\overline{v}=-\widetilde{f}(x,t),\\
 &\overline{v}(0,t)=0,\overline{v}(1,t)=-\widetilde{d}(t),\\
 &\overline{v}(x,0)=0.\notag
\end{align*}
Then for every $(x,s)\in [0,1]\times [0,t]$ there holds
\begin{align}
- \widetilde{v}(x,s)=&\overline{v}(x,s) \notag\\
           \leq & \max\Big\{\max\limits_{s\in[0,t]}-\widetilde{d}(s),0\Big\}\notag\\
                &+\frac{1}{\mu}2^{\frac{5p-8}{2p-4}}\max\limits_{(x,s)\in [0,1]\times[0,t]}|\widetilde{f}(x,s)|.\label{++17''}
\end{align}
\eqref{++13} follows from \eqref{++16} and \eqref{++17''}.

Finally, we conclude that
\begin{align*}
&\max\limits_{(x,s)\in[0,1]\times [0,t]} |v(x,s)| \\
=&\max\limits_{(x,s)\in[0,1]\times [0,t]} \big|\widetilde{v}(x,s)e^{-\frac{mx}{2\mu}}\big|\\
\leq &e^{\frac{|m|}{2\mu}}\max\limits_{(x,s)\in[0,1]\times [0,t]} |\widetilde{v}(x,s)|\\
\leq &e^{\frac{|m|}{2\mu}}\bigg(\max\limits_{s\in [0,t]}|\widetilde{d}(s)|+\frac{1}{\mu}2^{\frac{5p-8}{2p-4}}\max\limits_{(x,s)\in[0,1]\times [0,t]}|\widetilde{f}(x,s)|\bigg)\\
=&e^{\frac{|m|}{2\mu}}\bigg( \max\limits_{s\in [0,t]}\big|e^{\frac{m}{2\mu}}d(s)\big| \\
 &\;\;\;\;\;\;\;\; +\frac{1}{\mu}2^{\frac{5p-8}{2p-4}}\max\limits_{(x,s)\in[0,1]\times [0,t]}\big|f(x,s)e^{\frac{mx}{2\mu}}\big|\bigg)\\
\leq &e^{\frac{|m|}{\mu}}\bigg(\max\limits_{s\in [0,t]}|d(s)|+\frac{1}{\mu}2^{\frac{5p-8}{2p-4}}\max\limits_{(x,s)\in[0,1]\times [0,t]}|f(x,s)|\bigg).
\end{align*}
\end{IEEEproof}
\subsection{Proof of ISS properties for 1-$D$ transport PDE (Theorem~\ref{Theorem 6})}
\begin{IEEEproof}[Proof of (i) in Theorem~\ref{Theorem 6}]
By Theorem~\ref{Corollary 8} and Theorem~\ref{Theorem 9}, we have
\begin{align*}
 &|u(x,t)|
 = |w(x,t)+v(x,t)| \leq |w(x,t)|+|v(x,t)|\\
\leq & e^{\frac{|m|}{\mu}}\sqrt{2+\frac{m^2}{2\mu^2}} \|u_0\|_{H^1(0,1)} e^{-(\frac{m^2}{4\mu}+n+2\mu)t} \\
     &+e^{\frac{|m|}{\mu}}\bigg(\max\limits_{s\in [0,t]}|d(s)|+ \frac{1}{\mu}2^{\frac{5p-8}{2p-4}}\max\limits_{(x,s)\in[0,1]\times [0,t]}|f(x,s)|\bigg)\\
\leq & e^{\frac{|m|}{\mu}}\sqrt{2+\frac{m^2}{2\mu^2}} \|u_0\|_{H^1(0,1)} e^{-(\frac{m^2}{4\mu}+n+2\mu)t}\\
     &+e^{\frac{|m|}{\mu}}\bigg(\max\limits_{s\in [0,t]}|d(s)|+\frac{1}{\mu} 2^{\frac{5p-8}{2p-4}}\max\limits_{(x,s)\in[0,1]\times [0,t]}|f(x,s)|\bigg).
\end{align*}

\end{IEEEproof}

\begin{IEEEproof}[Proof of (ii) in Theorem~\ref{Theorem 6}]
Arguing as in the proof of Theorem~\ref{Corollary 8} and Theorem~\ref{Theorem 9}, letting $w(x,t)=\widetilde{w}(x,t)e^{-\frac{mx}{2\mu}}$ and $v(x,t)=\widetilde{v}(x,t)e^{-\frac{mx}{2\mu}}$, we consider the following two systems:
\begin{subequations}\label{++914}
 \begin{align}\label{++914.1}
 &\widetilde{w}_t-\mu \widetilde{w}_{xx}+\widetilde{n}\widetilde{w}=\widetilde{f}(x,t),\\
 &\widetilde{w}(0,t)=\widetilde{w}(1,t)=0,\label{++914.2}\\
  &\widetilde{w}(x,0)=\widetilde{w}_0,\label{++914.3}
\end{align}
\end{subequations}
and
\begin{subequations}\label{++12'}
 \begin{align}\label{++12.1'}
 &\widetilde{v}_t-\mu \widetilde{v}_{xx}+\widetilde{n}\widetilde{v}=0,\\
 &\widetilde{v}(0,t)=0,\widetilde{v}(1,t)=\widetilde{d}(t),\label{++12.2'}\\
  &\widetilde{v}(x,0)=0,\label{++12.3'}
\end{align}
\end{subequations}
where $\widetilde{n}=\frac{m^2}{4\mu}+n$, $\widetilde{w}_0=e^{\frac{mx}{2\mu}}u_0(x) $, $\widetilde{f}(x,t)=e^{\frac{mx}{2\mu}}f(x,t)$ and $\widetilde{d}(x,t)=e^{\frac{mx}{2\mu}}d(x,t)$.
Due to $u=w+v$, it suffices to estimate the solution of \eqref{++914} and \eqref{++12'}. Indeed, for \eqref{++12'}, we have by \eqref{++13}
\begin{align*}
\max\limits_{(x,s)\in[0,1]\times [0,t]} |\widetilde{v}(x,s)|\leq \max\limits_{s\in [0,t]}|\widetilde{d}(s)|.
\end{align*}

For \eqref{++914}, one may establish $H^1$-estimate of the solution as in the proof of Theorem~\ref{Theorem 7}, which implies its $L^\infty$-estimate by (i) of Lemma~\ref{Lemma 3}.

Indeed, multiplying \eqref{++914} by $\widetilde{w}$ and integrating over $(0,1)$, one may get
\begin{align*}
 \int_{0}^{1}\widetilde{w}_t\widetilde{w}\text{d}x+\mu\int_{0}^{1}\widetilde{w}_x^2\text{d}x+\widetilde{n}\int_{0}^{1}\widetilde{w}^2\text{d}x=
 \int_{0}^{1}\widetilde{f}\widetilde{w}\text{d}x,
\end{align*}
from which and Young's inequality, it yields for any $\varepsilon>0$,
\begin{align*}
 \frac{1}{2}\frac{d}{dt}\|\widetilde{w}\|^2+\mu\|\widetilde{w}_x\|^2+\widetilde{n}\|\widetilde{w}\|^2\leq \frac{1}{2\varepsilon}\|\widetilde{f}(\cdot,t)\|^2+\frac{\varepsilon}{2}\|\widetilde{w}\|^2.
\end{align*}
Note that $ \|\widetilde{w}_x\|^2\geq 2\|\widetilde{w}\|^2$. It follows
\begin{align*}
 \frac{d}{dt}\|\widetilde{w}(\cdot,t)\|^2+2(2\mu+\widetilde{n}- \frac{\varepsilon}{2})\|\widetilde{w}(\cdot,t)\|^2\leq  \frac{1}{\varepsilon}\|\widetilde{f}(\cdot,t)\|^2.
\end{align*}
By Gronwall's inequality, we may obtain
\begin{align}
     & \|\widetilde{w}(\cdot,t)\|^2 \notag\\
\leq & \|\widetilde{w}_0\|^2e^{-2(2\mu+\widetilde{n}- \frac{\varepsilon}{2})t}+\frac{1}{\varepsilon} \int_{0}^t\|\widetilde{f}(\cdot,s)\|^2\text{d}s \notag\\
\leq & \|e^{\frac{mx}{2\mu}}u_0\|^2e^{-2(2\mu+\widetilde{n}- \frac{\varepsilon}{2})t}+\frac{1}{\varepsilon}
           \int_{0}^t\big\|e^{\frac{m}{2\mu}}f(\cdot,s)\big\|^2\text{d}s\notag\\
\leq & e^{\frac{|m|}{\mu}}\|u_0\|^2e^{-2(2\mu+\widetilde{n}- \frac{\varepsilon}{2})t}+e^{\frac{|m|}{\mu}}\frac{1}{\varepsilon} \int_{0}^t\|f(\cdot,s)\|^2\text{d}s.\label{9151}
\end{align}

Multiplying \eqref{++914} by $\widetilde{w}_{xx}$ and integrating over $(0,1)$, one may get
\begin{align*}
 \int_{0}^{1}\widetilde{w}_{tx}\widetilde{w}_x\text{d}x+
 \mu\int_{0}^{1}\widetilde{w}_{xx}^2\text{d}x+
 \widetilde{n}\int_{0}^{1}\widetilde{w}_x^2\text{d}x
 =\int_{0}^{1}\widetilde{f}\widetilde{w}_{xx}\text{d}x,
\end{align*}
which yields for any $\varepsilon>0$,
\begin{align*}
 \frac{1}{2}\frac{d}{dt}\|\widetilde{w}_x\|^2
 +\mu\|\widetilde{w}_{xx}\|^2+\widetilde{n}\|\widetilde{w}_x\|^2\leq \frac{1}{2\varepsilon}\|\widetilde{f}(\cdot,t)\|^2+\frac{\varepsilon}{2}\|\widetilde{w}_{xx}\|^2.
\end{align*}
Arguing as in the proof of (i), it follows
\begin{align*}
 \frac{d}{dt}\|\widetilde{w}_x\|^2+
 2(2\mu+\widetilde{n}-\frac{\varepsilon}{2})\|\widetilde{w}_x\|^2\leq \frac{1}{2\varepsilon}\|\widetilde{f}(\cdot,t)\|^2.
\end{align*}
By Gronwall's inequality, we may obtain
\begin{align}
 \|\widetilde{w}_x(\cdot,t)\|^2
 \leq & \|\widetilde{w}_{0x}\|^2e^{-2(2\mu+\widetilde{n}-\frac{\varepsilon}{2})t}+\frac{1}{\varepsilon} \int_{0}^t\|\widetilde{f}(\cdot,s)\|^2\text{d}s\notag\\
 \leq & \Big\|\frac{m}{2\mu} e^{\frac{mx}{2\mu}} u_0+e^{\frac{mx}{2\mu}} u_{0x}\Big\|^2e^{-2(2\mu+\widetilde{n}-\frac{\varepsilon}{2})t}\notag\\
      & +\frac{1}{\varepsilon}\int_{0}^t\big\|e^{\frac{m}{2\mu}}f(\cdot,s)\big\|^2\text{d}s\notag\\
 \leq & e^{\frac{|m|}{\mu}}\bigg(\frac{m^2}{2\mu^2}  \| u_0\|^2+2\| u_{0x}\|^2\bigg)e^{-2(2\mu+\widetilde{n}-\frac{\varepsilon}{2})t} \notag\\
      &+e^{\frac{|m|}{\mu}}\frac{1}{\varepsilon} \int_{0}^t\|f(\cdot,s)\|^2\text{d}s.\label{9152}
\end{align}

\end{IEEEproof}

\begin{remark}
For the equation
\begin{align*}
u_t=\mu u_{xx}-mu_{x}-nu,
\end{align*}
where $\mu>0,m\geq 0, n\in\mathbb{R} $, with

Case 1) Dirichlet boundary conditions
\begin{align*}
&u(0,t)=d(t),\\
&u(1,t)=0,
\end{align*}
and

 Case 2) Robin (or Neumann) boundary conditions
\begin{align*}
&u(t,0)=d(t),\\
&u_x(1,t)=\left(\frac{m}{2\mu}-a\right)u(1,t),
\end{align*}
where $a\geq 0$, under the assumption $\frac{m^2}{4\mu}+n\geq 0$,
the ISS in $L^2$-norm is established in \cite{Karafyllis:2016} by the technique of eigenfunction expansions. Under the same assumption $\frac{m^2}{4\mu}+n\geq 0$, the ISS in $L^2$-norm in Case 2) is established in \cite{Zheng:2017} by the method of Lyapunov functionals. In the Case 1) with $m=0$, the ISS in $L^p$-norm ($p\in (2,+\infty)$) is established in \cite{Mironchenko:2017} by the maximum principles and monotonicity-based method. The ISS in $L^\infty$-norm (with $u(0,t)=d_0(t),u(1,t)=d_1(t)$) is established in \cite{karafyllis2017siam} by the technique of eigenfunction expansions.
\end{remark}

\begin{remark}
The method developed in this paper can be applied to linear problems with multidimensional spatial variables, e.g.,
\begin{subequations}
 \begin{align*}
 &u_t-\mu \Delta u+c(x,t)u=f(x,t),\ \ \text{in }\  \Omega\times \mathbb{R}_{\geq 0},\\
 &u(x,t)=d(t),\ \ \text{in }\ \partial\Omega\times \mathbb{R}_{\geq 0},\label{}\\
  &u(x,0)=u_0(x),\ \ \text{in }\  \Omega,
\end{align*}
\end{subequations}
where $\Omega\subset\mathbb{R}^n (n\geq 2)$ is an open bounded domain with smooth boundary $\partial \Omega$, $ c(x,t)$ is a smooth function in $ \Omega\times \mathbb{R}_{\geq 0}$ with $0< m\leq c(x,t)\leq M$, $\Delta$ is the Laplace operator, $\mu>0$ is a constant, and $u_0\in H^2(\Omega)$. Due to
\begin{equation*}
H^1(\Omega)\hookrightarrow \left\{
  \begin{array}{ll}\nonumber
  L^{q}(\Omega),\ \ \forall q\in [1,+\infty),\ n=2; \\
    L^{q}(\Omega),\ \ \forall q\in [1,\frac{2n}{n-2}],\ n>2;
  \end{array}\notag
\right.
\end{equation*}
one may get
\begin{align*}
\|u(\cdot,t)\|_{L^q(\Omega)}\leq & C_0\|u_0\|_{H^1(\Omega)} e^{-\lambda t}
+ C_1\max\limits_{s\in [0,t]}|d(s)| \notag\\
&+ C_2 \max\limits_{(x,s)\in\overline{Q}_t}|f(x,s)|,
\end{align*}
where $C_0,C_1,C_2,\lambda$ are universal positive constants.

 Indeed, one may obtain the $H^1$-estimate (bounded by $C\|u_0\|_{H^1(\Omega)} e^{-\lambda t}$) of the solution $w$ of the following system
 \begin{align*}
 &w_t-\mu \Delta w+c(x,t)w=0,\ \ \text{in }\  \Omega\times \mathbb{R}_{\geq 0},\\
 &w(x,t)=0,\ \ \text{in }\  \partial\Omega\times \mathbb{R}_{\geq 0},\label{}\\
  &w(x,0)=u_0(x),\ \ \text{in }\  \Omega.
\end{align*}
The $L^{q}$-estimate of $w$ follows from the embedding theorem. Then by De~Giorgi iteration, one may get the maximum estimate (bounded by $C\max\limits_{s\in [0,t]}|d(s)|+ C \max\limits_{(x,s)\in\overline{Q}_t}|f(x,s)|$) of $v$, which is the solution of the following system
\begin{align*}
 &v_t-\mu \Delta v+c(x,t)v=f(x,t),\ \ \text{in }\  \Omega\times \mathbb{R}_{\geq 0},\\
 &v(x,t)=d(t),\ \ \text{in }\  \partial\Omega\times \mathbb{R}_{\geq 0},\label{}\\
  &v(x,0)=0,\ \ \text{in }\  \Omega.
\end{align*}
Due to $L^\infty(\Omega)\hookrightarrow L^q(\Omega)$, one may get the $L^{q}$-estimate of $v$. Finally, the $L^{q}$-estimate of $u$ is given by the $L^{q}$-estimates of $w$ and $v$.

%

\end{remark}
\subsection{ISS of a 1-$D$ transport equation with anti-stable term under boundary feedback control}
As an application of Theorem~\ref{Theorem 6}, we study the stability of an anti-stable 1-$D$ transport equation:
\begin{align}\label{++9181}
 u_t-\mu u_{xx}+a(x)u=f(x,t) \ \ \ \ \text{in}\ \  (0,1)\times \mathbb{R}_{\geq 0},
\end{align}
where $\mu >0$ is a constant, $a\in C^1([0,1])$ and $f(x,t)\in C^1([0,1]\times \mathbb{R}_{\geq 0})$, subject to the boundary and initial conditions
\begin{subequations}\label{++9182}
\begin{align}
 &u(0,t)=0,u(1,t)=U(t),\\
 &u(x,0)=u_0(x),
\end{align}
\end{subequations}
where $U(t)\in \mathbb{R}$ is the control input.

The stabilization of \eqref{++9181} in a disturbance-free setting with $\mu =1$ and $f(x,t)\equiv 0$ is presented in \cite{Krstic:2008,Liu:2003,Smyshlyaev:2004}, and the ISS properties w.r.t. boundary disturbances have been addressed in \cite{Karafyllis:2016,Karafyllis:2016a,Mironchenko:2017}. It should be noticed that in different work, the control input is placed on the left or right end of the system. Nevertheless, it can be switched to the other end by a spatial variable transformation $x\rightarrow 1-x$.

In \cite{Krstic:2008,Liu:2003,Smyshlyaev:2004}, the exponential stability of parabolic PDEs of the form \eqref{++9181}, with $\mu =1$ and $f(x,t)\equiv 0$, has been achieved by means of a boundary feedback control of the form
\begin{align}
U(t)=\int_{0}^1k(1,y)u(y,t)\text{d}y,\ \forall t\geq 0\label{++9183}
\end{align}
where $k\in C^2([0,1]\times[0,1])$ is an appropriate function.

In fact, as $\mu >0$ and $f(x,t)\in C^1([0,1]\times \mathbb{R}_{\geq 0})$, the function $k$ can also be obtained as the Volterra kernel of a Volterra integral transformation
\begin{align*}
w(x,t)=u(x,t)-\int_{0}^xk(x,y)u(y,t)\text{d}y,\ \forall (x,t)\in[0,1]\times \mathbb{R},
\end{align*}
which transforms \eqref{++9181}, \eqref{++9182}, and \eqref{++9183} to the problem \begin{align}\label{++9185}
 w_t-\mu w_{xx}+nw=f(x,t) \ \ \ \ \text{in}\ \  (0,1)\times \mathbb{R}_{\geq 0},
\end{align}
with $n\geq 0$, subject to the boundary and initial conditions
\begin{align*}
&w(0,t)=w(1,t)=0,\\
&w(x,0)=w_0(x)=u_0(x)-\int_{0}^xk(x,y)u_0(y)\text{d}y.
\end{align*}
The free parameter $n\geq 0$ can be used to
set the convergence rate. The solution of the original problem can be found by the
inverse Volterra integral transformation
\begin{align}
u(x,t)=w(x,t)+\int_{0}^xl(x,y)w(y,t)\text{d}y,\ \forall (x,t)\in[0,1]\times \mathbb{R},\label{++9187}
\end{align}
where $l\in C^2([0,1]\times[0,1])$ is an appropriate kernel. The existence of
the kernels $k\in C^2([0,1]\times[0,1])$ and $l\in C^2([0,1]\times[0,1])$ can be obtained in the same way as in \cite{Liu:2003,Smyshlyaev:2004}.

 In the presence of actuator errors, i. e., when the applied control action is
of the form
\begin{align}
U(t)=d(t)+\int_{0}^1k(1,y)u(y,t)\text{d}y,\ \forall t\geq 0,\label{++9188}
\end{align}
where $d\in C^2 (\mathbb{R}_{\geq 0}) $, the transformed solution $ w(x,t)$ satisfies \eqref{++9185} subject to the
boundary and initial conditions
\begin{subequations}\label{++9189}
\begin{align}
&w(0,t)=0,w(1,t)=d(t),\\
&w(x,0)=w_0(x)=u_0(x)-\int_{0}^xk(x,y)u_0(y)\text{d}y.
\end{align}
\end{subequations}
According to Theorem~\ref{Theorem 6}, we have the following ISS estimates for \eqref{++9185} with \eqref{++9189}:
\begin{align*}
\max_{x\in[0,1]}|w(x,t)|\leq & \sqrt{2}\|w_0\|_{H^1(0,1)}e^{-(n+2\mu) t}+\max\limits_{s\in [0,t]}|d(s)|\notag\\
   &+2^{\frac{5p-8}{2p-4}}\frac{1 }{\mu}\max\limits_{(x,s)\in[0,1]\times [0,t]}|f(x,s)|,\ \forall p>2,
\end{align*}
and
\begin{align*}
\max_{x\in[0,1]}|w(x,t)| \leq &\sqrt{2}\|w_0\|_{H^1(0,1)}e^{-(n+2\mu-\frac{\varepsilon}{2}) t}+\max\limits_{s\in [0,t]}|d(s)|\notag\\
 &+\sqrt{\frac{3}{\varepsilon} \int_{0}^t\|f(\cdot,s)\|^2\text{d}s},\ \forall \varepsilon \in (0, 2n+4\mu).
\end{align*}

Note that
\begin{align*}
 \|w_0\|_{H^1(0,1)}&\leq \|u_0-\int_{0}^xk(x,y)u_0(y)\text{d}y \|_{H^1(0,1)}\notag\\
 &\leq \|u_0\|_{H^1(0,1)}+\|\int_{0}^xk(x,y)u_0(y)\text{d}y \|_{H^1(0,1)}\notag\\
 &\leq \|u_0\|_{H^1(0,1)}+2\|k(\cdot,\cdot)u_0(\cdot) \|\notag\\
 &\leq C\|u_0\|_{H^1(0,1)},
\end{align*}
where $C$ depends only on $ \max\limits_{x\in[0,1] }|k(x,x)|$.

Finally, by \eqref{++9187}, we have the following ISS estimates for systems \eqref{++9181} and \eqref{++9182} w.r.t. control actuator errors for boundary
state feedback \eqref{++9188}:
\begin{align*}
 \max_{x\in[0,1]}|u(x,t)|\leq & \bigg(1+\max_{(x,y)\in[0,1]\times[0,1] }|l(x,y)|\bigg)\max_{x\in[0,1]}|w(x,t)|\\
 \leq & C\|u_0\|_{H^1(0,1)} e^{-(n+2\mu) t}+C\max\limits_{s\in [0,t]}|d(s)| \\
 &+\frac{C }{\mu}2^{\frac{5p-8}{2p-4}}\max\limits_{(x,s)\in[0,1]\times [0,t]}|f(x,s)|,\ \ \forall p>2,
\end{align*}
and
\begin{align*}
 \max_{x\in[0,1]}|u(x,t)|\leq & C\|u_0\|_{H^1(0,1)}e^{-(n+2\mu-\frac{\varepsilon}{2}) t}+C\max\limits_{s\in [0,t]}|d(s)|\notag\\
 &+C\sqrt{\frac{3}{\varepsilon} \int_{0}^t\|f(\cdot,s)\|^2\text{d}s}, \forall \varepsilon \in (0, 2n+4\mu),
\end{align*}
where $C$ is a constant depending only on $ \max\limits_{x\in[0,1] }|k(x,x)|$ and $ \max\limits_{(x,y)\in[0,1]\times[0,1] }|l(x,y)|$.

\begin{remark} In the case where $a(x)\equiv a$ is a constant, the ISS in $L^2$-norm and $L^p$-norm ($p>2$) for System \eqref{++9181} and \eqref{++9182} w.r.t. control actuator errors for boundary feedback control \eqref{++9188} is established in \cite{Karafyllis:2016a} and \cite{Mironchenko:2017}, respectively. Moreover, a continuous controller $U(t)$ is designed to stabilize the system \eqref{++9181} and \eqref{++9182} in $L^2$-state space with disturbance $d(t)$ on the boundary under the assumption that $ d(t)$ and its derivative $d'(t)$ are uniformly bounded on $ \mathbb{R}_{\geq 0}$.
\end{remark}
\begin{remark} Consider the system
\begin{align}\label{++9191}
 u_t-\mu u_{xx}+mu_x+nu=f(x,t) \ \ \ \ \text{in}\ \  (0,1)\times \mathbb{R}_{\geq 0},
\end{align}
where $\mu >0$ and $m,n\in\mathbb{R}$ are constants, $f(x,t)\in C^1([0,1]\times \mathbb{R}_{\geq 0})$, subject to the boundary and initial conditions
\begin{subequations}\label{++9192}
\begin{align}
 &u(0,t)=0,u(1,t)=U(t),\\
 &u(x,0)=u_0(x),
\end{align}
\end{subequations}
where $U(t)\in \mathbb{R}$ is the control input. If we design the feedback as
\begin{align}
U(t)=e^{-\frac{m}{2\mu}}\bigg(d(t)+\int_{0}^1k(1,y)e^{\frac{my}{2\mu}}u(y,t)\text{d}y\bigg),\ \forall t\geq 0,\label{++9193}
\end{align}
where $d\in C^2 (\mathbb{R}_{\geq 0}) $, one may obtain the ISS properties as well. Indeed, letting $u(x,t)=e^{-\frac{mx}{2\mu}}v(y,t) $, we get
\begin{subequations}\label{++9194}
\begin{align}
&v_t-\mu v_{xx}+\widetilde{n}v=\widetilde{f}(x,t),\\
 &v(0,t)=0,v(1,t)=d(t)+\int_{0}^1k(1,y)v(y,t)\text{d}y,\\
 &v(x,0)=e^{\frac{mx}{2\mu}}u_0(x),
\end{align}
\end{subequations}
where $\widetilde{n}=\frac{m^2}{4\mu}+n$  and $ \widetilde{f}(x,t)=e^{\frac{mx}{2\mu}}f(x,t)$.
The ISS properties for \eqref{++9194} can be established as \eqref{++9181}, \eqref{++9182}, and \eqref{++9188}. Then the ISS properties for \eqref{++9191}, \eqref{++9192}, and \eqref{++9193} can be obtained by the transformation $u(x,t)=e^{-\frac{mx}{2\mu}}v(y,t) $.
\end{remark}

\section{ISS Properties for Burgers' Equation}\label{Sec: Burgers Eq}

%
%
%
In this section, we establish the ISS properties for Burgers' equation w.r.t. boundary and in-domain disturbances described in Theorem~\ref{Theorem 11}.

To split the original problem, we consider the following two systems:
\begin{subequations}\label{++29}
 \begin{align}
 &w_t-\mu w_{xx}+\nu ww_x=0\ \ \text{in}\ (0,1)\times \mathbb{R}_{\geq 0},\\
 &w(0,t)=0,w(1,t)=d(t),\\
  &w(x,0)=0,
\end{align}
\end{subequations}
and
\begin{subequations}\label{++31}
 \begin{align}
 &v_t-\mu v_{xx}+\nu vv_x+\nu wv_x +\nu vw_x=f(x,t), \\
 &\ \ \ \ \ \ \ \ \ \ \ \text{in}\ (0,1)\times \mathbb{R}_{\geq 0}, \nonumber\\
 &v(0,t)=v(1,t)=0,\\
 &v(x,0)=u_0(x).
\end{align}
\end{subequations}

For system \eqref{++29}, we have the following estimates.
\begin{theorem} \label{Theorem 12} For any $t>0$, there holds
\begin{align}\label{++32''}
\max\limits_{(x,s)\in[0,1]\times [0,t]} |w(x,s)|\leq \max\limits_{s\in [0,t]}|d(s)|.
\end{align}
It follows for all $t>0$,
\begin{align*}
\|w(\cdot,t)\|\leq \max\limits_{s\in [0,t]}|d(s)|.
\end{align*}
\end{theorem}
For system \eqref{++31}, we have the following estimate.
\begin{theorem} \label{Theorem 13} Suppose that $\nu>0$ and $\max\limits_{t\in \mathbb{R}_{\geq 0}} |d(t)|<\frac{\mu}{\nu}$. For every $t>0$, there holds
\begin{align*}
\|v(\cdot,t)\|^2\leq \|u_0\|^2 e^{-(\mu -\varepsilon)t}+\frac{1}{\varepsilon} \int_{0}^t\|f(\cdot,s)\|^2\text{d}s,\ \forall \varepsilon\in (0,\mu).
\end{align*}
\end{theorem}
\begin{IEEEproof}[Proof of Theorem~\ref{Theorem 12}]
One may proceed as the proof of Theorem~\ref{Theorem 9}. For any fixed $t>0$, let $k_0,k$, $ \eta(x,s)$ and $t_0$ be defined as in the proof of Theorem~\ref{Theorem 9}.
Taking $\eta(x,s)$ as a test function and proceeding as before, it suffices to estimate the integration $ \int_{0}^1w(x,t_0)w_x(x,t_0)(w(x,t_0)-k)_+ \text{d}x$. We write $w=w(x,t_0)$ for simplicity.

One may get
\begin{align*}
&\int_{0}^1ww_x(w-k)_+ \text{d}x \\
=& \int_{0}^1(w-k)_+((w-k)_+)_x(w-k)_+ \text{d}x \\
 &+\int_{0}^1k((w-k)_+)_x(w-k)_+ \text{d}x\\
=& \frac{1}{3}(w-k)_+^{3}|^{x=1}_{x=0}+\frac{k}{2}(w-k)_+^{2}|^{x=1}_{x=0}=0.
\end{align*}
Finally, one may obtain \eqref{++32''} by repeating the process of the proof of Theorem~\ref{Theorem 9}.
\end{IEEEproof}
\begin{IEEEproof}[Proof of Theorem~\ref{Theorem 13}]
Multiplying \eqref{++31} by $v$ and integrating over $(0,1)$, we get
\begin{align*}
&\int_{0}^1v_tv\text{d}x+\mu\int_{0}^1v^2_{x}\text{d}x+\nu\int_{0}^1v^2v_x\text{d}x+\nu\int_{0}^1(wv)_xv\text{d}x \\
=&\int_{0}^1f(x,t)v\text{d}x.
\end{align*}
Note that $\int_{0}^1v^2v_x\text{d}x=\frac{1}{3}v^3|^{x=1}_{x=0}=0$ and
\begin{align*}
\int_{0}^1(wv)_xv\text{d}x = wv^2 |^{x=1}_{x=0}-\int_{0}^1wvv_x\text{d}x=-\int_{0}^1wvv_x\text{d}x.
\end{align*}
By Young's inequality, H\"{o}lder inequality, Theorem~\ref{Theorem 12}, and the assumption on $d$, we deduce that
\begin{align*}
&\frac{1}{2}\frac{d}{dt}\|v\|^2+\mu\|v_x\|^2 \\
\leq &\nu\int_{0}^1|wvv_x|\text{d}x+\int_{0}^1f(x,t)v\text{d}x\\
\leq &\frac{\nu}{2}\max\limits_{(x,s)\in[0,1]\times [0,t]} |w(x,s)|(\|v\|^2+\|v_x\|^2) \\
     &+\frac{1}{2}(\|f(\cdot,t)\|^2+\|v\|^2)\\
\leq &\frac{\nu}{2}\max\limits_{s\in [0,t]}|d(s)|(\|v\|^2+\|v_x\|^2)+\frac{1}{2\varepsilon}\|f(\cdot,t)\|^2
+\frac{\varepsilon}{2}\|v\|^2\\
\leq & \frac{\nu}{2}\frac{\mu}{\nu}(\|v\|^2+\|v_x\|^2)
+\frac{1}{2\varepsilon}\|f(\cdot,t)\|^2+\frac{\varepsilon}{2}\|v\|^2\\
=&\frac{1}{2}(\varepsilon+\mu)\|v\|^2
+\frac{\mu}{2}\|v_x\|^2+\frac{1}{2\varepsilon}\|f(\cdot,t)\|^2,
\end{align*}
where we may choose $0<\varepsilon <\mu$.

By Lemma~\ref{Lemma 2}, we have
\begin{align*}
\mu\|v_x\|^2=\frac{\mu}{2}\|v_x\|^2+\frac{\mu}{2}\|v_x\|^2\geq \frac{\mu}{2}\|v_x\|^2+\mu\|v\|^2.
\end{align*}

Then we have
\begin{align*}
\frac{d}{dt}\|v(\cdot,t)\|^2
\leq-(\mu -\varepsilon)\|v(\cdot,t)\|^2+\frac{1}{\varepsilon}\|f(\cdot,t)\|^2.
\end{align*}
By Growall's inequality, we have
\begin{align*}
\|v(\cdot,t)\|^2
&\leq\|v(\cdot,0)\|^2 e^{-(\mu -\varepsilon)t}+\frac{1}{\varepsilon} \int_{0}^t\|f(\cdot,s)\|^2\text{d}s\\
&=\|u_0\|^2 e^{-(\mu -\varepsilon)t}+\frac{1}{\varepsilon} \int_{0}^t\|f(\cdot,s)\|^2\text{d}s.
\end{align*}
\end{IEEEproof}
\begin{IEEEproof}[Proof of (i) in Theorem~\ref{Theorem  11}] Note that $u=w+v$, we get by Theorem~\ref{Theorem  12} and Theorem~\ref{Theorem  13}:
\begin{align*}
\|u(\cdot,t)\|^2
\leq & 2\|w(\cdot,t)\|^2+2\|v(\cdot,t)\|^2\\
\leq & 2\|u_0\|^2 e^{-(\mu -\varepsilon)t}+2(\max\limits_{s\in [0,t]}|d(s)|)^2 \\
     &+\frac{2}{\varepsilon} \int_{0}^t\|f(\cdot,s)\|^2\text{d}s,\ \forall \varepsilon\in (0,\mu).
\end{align*}
\end{IEEEproof}
In order to prove (ii) in Theorem~\ref{Theorem 11}, one may consider the following two systems.
\begin{subequations}\label{+++32}
 \begin{align}
 &w_t-\mu w_{xx}+\nu ww_x=f(x,t)\ \ \text{in}\  (0,1)\times \mathbb{R}_{\geq 0},\\
 &w(0,t)=0,w(1,t)=d(t),\\
  &w(x,0)=0,
\end{align}
\end{subequations}
and
\begin{subequations}\label{+++33}
 \begin{align}
 &v_t-\mu v_{xx}+\nu vv_x+\nu wv_x +\nu vw_x=0\ \ \text{in}\ (0,1)\times \mathbb{R}_{\geq 0},\\
 &v(0,t)=v(1,t)=0,\\
  &v(x,0)=u_0(x),
\end{align}
\end{subequations}

For system \eqref{+++32}, we have the following estimates.
\begin{theorem} \label{Theorem 15}For every $t>0$, there holds for any $p\in(2,+\infty)$,
\begin{align*}
&\max\limits_{(x,s)\in[0,1]\times [0,t]} |w(x,s)| \nonumber\\
\leq &\max\limits_{s\in [0,t]}|d(s)|+\frac{1}{\mu}2^{\frac{5p-8}{2p-4}}\max\limits_{(x,s)\in[0,1]\times [0,t]}|f(x,s)|.
\end{align*}
It follows that for any $p\in(2,+\infty)$,
\begin{align*}
     &\|w(\cdot,t)\| \nonumber\\
\leq & \max\limits_{s\in [0,t]}|d(s)|+\frac{1}{\mu}2^{\frac{5p-8}{2p-4}}\max\limits_{(x,s)\in[0,1]\times [0,t]}|f(x,s)|.
\end{align*}
\end{theorem}

For system \eqref{+++33}, we have the following estimate.
\begin{theorem} \label{Theorem 16}For every $t>0$, there holds
\begin{align*}
\|v(\cdot,t)\|^2\leq \|u_0\|^2 e^{-\mu t}.
\end{align*}
\end{theorem}
Then the result of (ii) in Theorem~\ref{Theorem 11} is a consequence of Theorem~\ref{Theorem 15} and Theorem~\ref{Theorem 16} under the assumption that $\max\limits_{t\in \mathbb{R}_{\geq 0}} |d(t)|+ \frac{1}{\mu}2^{\frac{5p-8}{2p-4}}\max\limits_{(x,s)\in[0,1]\times \mathbb{R}_{\geq 0}}|f(x,s)|<\frac{\mu}{\nu}$ for some $p\in (2,+\infty)$. As the development can proceed in a similar way as above, the details on the proof of Theorem~\ref{Theorem 15}, Theorem~\ref{Theorem 16}, and (ii) of Theorem~\ref{Theorem 11} are omitted.
\begin{remark}
In general, the boundness of the disturbances is a reasonable assumption for nonlinear PDEs in the establishment of ISS properties \cite{Mironchenko:2016}.
\end{remark}
\section{Concluding Remarks}\label{Sec: Conclusion}
This paper applied the technique of De~Giorgi iteration to the establishment of ISS properties w.r.t. boundary and in-domain disturbances for some semi-linear PDEs with Dirichlet boundary conditions. The ISS in $L^\infty$-norm for 1-$D$ transport PDEs and the ISS in $L^2$-norm for Burgers' equation have been obtained. As pointed out in Section~\ref{Sec: Transport PDE}, the method developed in this paper can be generalized to some problems on multidimensional spatial domain. Moreover, some estimates established in this paper are in $L^\infty$-norm, which are the enhancements of $L^p$-norm ($p\in(2,+\infty)$) obtained in \cite{Mironchenko:2017}. It is also a complement of \cite{Zheng:2017} where the ISS in $L^2$-norm has been established for some 1-$D$ transport PDEs with Robin (or Neumann) boundary conditions. As the De~Giorgi iteration is a well-established tool for regularity analysis of PDEs, we can expect that the method developed in this work can be extended to the study of a wider class of nonlinear parabolic PDEs. Future work may include the study of Burgers' equation with boundary feedback control. The investigation on ISS properties for other well-known semi-linear parabolic equations, such as generalized Burgers' equations, Fisher-Kolmogorov equation, and Chaffee-Infante equation, may also be of theoretical and practical interest.



%
%

\ifCLASSOPTIONcaptionsoff
  \newpage
\fi



\bibliographystyle{IEEEannot}
\bibliography{IEEEabrv,References}

\vspace*{-2\baselineskip}

\end{document}